\documentclass{etna}
\usepackage{graphicx}
\usepackage{amsmath}
\usepackage{amsfonts}
\usepackage{MnSymbol}

\def\smallskip{\vskip 5pt plus1pt minus1pt}
\def\medskip{\vskip 9pt plus2pt minus2pt}
\def\bigskip{\vskip 18pt plus4pt minus4pt}

\def\sqr#1#2{{\vcenter{\vbox{\hrule height.#2pt
        \hbox{\vrule width.#2pt height#1pt \kern#1pt \vrule width.#2pt}
             \hrule height.#2pt}}}}

\usepackage{lineno}
\begin{document}

\thispagestyle{empty}
\bibliographystyle{siam}
\title{A note on parallel preconditioning for all-at-once evolutionary PDEs}

\author{Anthony Goddard\thanks{Mathematical Institute, Oxford University, UK. Supported 
by the James Pantyfedwen Foundation.} \and
Andy Wathen\thanks{Mathematical Institute, Oxford University, UK}}

\date{}

\maketitle

\markboth{A.~J.~Goddard and A.~J.~Wathen}{Parallel preconditioning for
evolutionary PDEs}

\begin{abstract}
McDonald, Pestana and 
Wathen (SIAM J. Sci. Comput. {\bf 40}(2), pp.~A2012--A1033, 2018) 
present a method for
preconditioning of time-dependent PDEs via approximation 
by a nearby time-periodic problem, that is, they employ
circulant-related matrices as preconditioners for
the non-symmetric block Toeplitz matrices which arise
from an all-at-once formulation. They suggest that
such an approach might be efficiently implemented
in parallel. In this short article, we present parallel
numerical results for their preconditioner which 
exhibit strong scaling. We also extend their 
preconditioner via a Neumann series approach, which also
allows for efficient parallel execution.   

Our simple implementation (in C++ and MPI) is available 
at the Git repository \texttt{PARALAAOMPI}.\footnote{\texttt{https://github.com/anthonyjamesgoddard/PARALAAOMPI}}
\end{abstract}


\begin{AMS} 65M20,65F08,65Y05.
\end{AMS}

\section{Introduction}\label{intro}
There have been several suggestions for ways to achieve parallel
computationally efficient methods for time-dependent problems.
Perhaps the most common approaches are based on the Parareal 
algorithm \cite{MADAY} and its use together with multilevel
ideas \cite{FALG}, although there are several
other approaches; see the review by Gander \cite{GAND}. A recent suggestion by
McDonald et al. \cite{McDonald_Pestana_W18} involves the use
of circulant-related preconditioners for the block Toeplitz matrices
that arise from the approximation of initial value problems with constant
time-steps. Their approach is particularly geared to linear initial value
problems for PDEs, although it can be applied in the simpler context of 
ODE IVPs \cite{HON}. For PDEs, regularity of the spatial grid is not required.

The approach requires the solution of block diagonal
 systems in
different orderings so as to effect the action of the
preconditioner, as we show below. It also allows for an extension
involving a Neumann series which we introduce here. Both the extended 
and the original preconditioners would appear suitable for 
effective parallel implementation, although such implementation
details have not been explored before. That the original preconditioner
leads to a small number of iterations which is independent of the number
of time-steps when employed with 
the widely used GMRES method \cite{Saad_Schultz86} is established in
\cite{McDonald_Pestana_W18}.  

In this short article we make a preliminary exploration of
the possibilities for effective parallel implementation of
these preconditioners. Our initial results---using C++ and MPI---show 
strong scaling with up to 32 cores for the preconditioned 
GMRES solution of the all-at-once (monolithic) system. This system is derived 
from a spatial finite element approximation and a simple time-stepping strategy in the usual
method of lines approach. We explore all-at-once formulations for the
heat equation and for the wave equation together with associated initial and boundary conditions.

In Section 2 we describe the original McDonald-Pestana-Wathen
preconditioner and our extension of it. The details of our
parallel implementation are given in Section 3 and numerical
(timing) results in Section 4, followed by conclusions.

\section{Description of the Preconditioners}
\label{preconditioners}
Consider the heat equation
\begin{linenomath*}
\begin{align}
\frac{\partial u}{\partial t} &= \frac{\partial^2 u}{\partial x^2} \quad \text{ on } \Omega \times (0,T],\label{heatbegin} \\
u &= 0 \qquad \text{ for } \mathbf{x} \in \partial \Omega, \\
u(0,\mathbf{x}) &= s(\mathbf{x}),
\label{heatend}
\end{align}
\end{linenomath*}

where $\Omega = [0,1]$ and $T = 1$. Discretising in space with standard Galerkin finite elements we obtain
\begin{linenomath*}
\begin{equation}
M\frac{\text{d} \mathbf{u}}{\text{d} t} = -K \mathbf{u},
\label{applyEuler}
\end{equation}
\end{linenomath*}

where $M,K \in \mathbb{R}^{n\times n}$ are the mass and stiffness matrices and $n$ is the number of nodes in the spatial discretisation. 
Now we use the implicit Euler scheme to discretise the temporal domain to obtain
\begin{linenomath*}
\begin{equation}
(K+\tau M) \mathbf{u}_k = M\mathbf{u}_{k-1} , \qquad k \in [1,\ell],
\label{timestepproblem1}
\end{equation}
\end{linenomath*}

where $\tau$ is the constant time-step, $\ell$ is the number of
time-steps in the temporal discretisation and $\mathbf{u}_0$ is a projection of
the initial data. The idea presented in \cite{McDonald_Pestana_W18} involves packaging the approximate solutions into a so-called monolithic system. Executing this idea yields
\begin{linenomath*}
\begin{equation}
\mathcal{A} \mathbf{U}=
\begin{bmatrix}
    A_0 &  &     &  \\
    A_1 & A_0 &   &  \\
     & \ddots & \ddots  &  \\
     &  &   A_1  & A_0
\end{bmatrix}
\begin{bmatrix}
    \mathbf{u}_1  \\
    \mathbf{u}_2   \\
    \vdots  \\
    \mathbf{u}_{\ell}
\end{bmatrix}
=
\begin{bmatrix}
    M\mathbf{u}_0 \\
    \mathbf{0}   \\
    \vdots  \\
    \mathbf{0}
\end{bmatrix} = \mathbf{b},
\label{monolithicheat}
\end{equation}
\end{linenomath*}

where $A_0 = K+\tau M$, $A_1=-M$ and $\mathcal{A} \in \mathbb{R}^{n\ell \times n\ell}$. 
We note that the monolithic system does not need to be formed explicitly; it is merely used as a conduit for demonstration purposes. The McDonald-Pestana-Wathen preconditioner in its original form is the block circulant matrix
\begin{linenomath*}
\begin{equation}
\mathcal{P}=
\begin{bmatrix}
    A_0 &  &     & A_1 \\
    A_1 & A_0 &   &  \\
     & \ddots & \ddots  &  \\
     &  &   A_1  & A_0
\end{bmatrix}.
\label{standardheatprecon}
\end{equation}
\end{linenomath*}

We can extend the all-at-once method to non-uniform time-stepping schemes. Consider applying the implicit Euler scheme to (\ref{applyEuler}) with variable time-steps $\tau_1,\tau_2,...,\tau_{\ell}$: Equation (\ref{timestepproblem1}) then becomes
\begin{linenomath*}
\begin{equation}
(K+\tau_i M) \mathbf{u}_k = M\mathbf{u}_{k-1}  , \qquad k \in [1,\ell].
\label{nonuniformtimestepproblem1}
\end{equation}
\end{linenomath*}

We can package this sequence into a monolithic system as follows:
\begin{linenomath*}
\begin{equation}
\mathcal{B} \mathbf{U}=
\begin{bmatrix}
    A_0^1 &  &     &  \\
    A_1 & A_0^2  &   &  \\
     & \ddots & \ddots  &  \\
     &  & A_1  & A_0^\ell
\end{bmatrix}
\begin{bmatrix}
    \mathbf{u}_1  \\
    \mathbf{u}_2   \\
    \vdots  \\
    \mathbf{u}_\ell
\end{bmatrix}
=
\begin{bmatrix}
    M\mathbf{u}_0 \\
    \mathbf{0}   \\
    \vdots  \\
    \mathbf{0}
\end{bmatrix} = \mathbf{b},
\label{nonuniformmonolithicheat}
\end{equation}
\end{linenomath*}

where $A_0^i = K+\tau_i M$ and $A_1=-M$. While we have lost the block Toeplitz structure of the system, this will not prevent us from applying
\begin{linenomath*}
\begin{equation}
\mathcal{Q}=
\begin{bmatrix}
    A_0^1 &  &     & A_1 \\
    A_1 & A_0^2 &    &  \\
     & \ddots & \ddots  &  \\
     &  &   A_1  & A_0^{\ell}
\end{bmatrix}
\label{nonunistandardheatprecon}
\end{equation}
\end{linenomath*}

as a preconditoner. The issue is that of the computational cost of the application of the preconditioner. We can write
\begin{linenomath*}
\begin{equation}
\mathcal{Q} = \mathcal{P} + \sigma \otimes K,
\end{equation}
\end{linenomath*}

where $\sigma$ is a diagonal matrix with diagonal entries given by
$\sigma_i = \tau_i - \tau$, $i\in [0,\ell ]$, and $\tau = 1/\ell$. Assuming that $||\sigma \otimes K|| < 1$ we can use a Neumann approximation to calculate the inverse of $\mathcal{Q}$ to obtain
\begin{linenomath*}
\begin{equation}
\mathcal{Q}_i^{-1} = \mathcal{P}^{-1} - \mathcal{P}^{-1}([\sigma \otimes K] \mathcal{P}^{-1}) + \mathcal{P}^{-1}([\sigma \otimes K]\mathcal{P}^{-1})^2 + ... (-1)^{i-1} \mathcal{P}^{-1}([\sigma \otimes K] \mathcal{P}^{-1})^{i-1},
\label{preconapprox1}
\end{equation}
\end{linenomath*}

where $i$ is some positive integer that we must choose.

We can also apply the all-at-once method to hyperbolic equations. We will restrict our scope to the wave equation
\begin{linenomath*}
\begin{align}
\frac{\partial^2 u}{\partial t^2} &= \frac{\partial^2 u}{\partial x^2}, \quad \text{ on } \Omega \times (0,T],
\label{wavebegin} \\
u &= 0, \qquad \text{ for } \mathbf{x} \in \partial \Omega, \\
\frac{\partial u}{\partial t}(0,\mathbf{x}) &= 0, \quad \\
u(0,\mathbf{x}) &= s(\mathbf{x}).
\label{waveend}
\end{align}
\end{linenomath*}

Discretising in space we obtain
\begin{linenomath*}
\begin{equation}
M\frac{\text{d}^2 \mathbf{u}}{\text{d} t^2} = -K \mathbf{u}.
\label{applyEuler2}
\end{equation}
\end{linenomath*}

We can choose the central difference formula to approximate the second time derivative, resulting in the sequence
\begin{linenomath*}
\begin{equation}
M\mathbf{u}_{n-1} + (\tau^2K-2M)\mathbf{u}_{n} + M\mathbf{u}_{n+1} = \mathbf{0}.
\label{CenDiffDisc}
\end{equation}
\end{linenomath*}

Casting this sequence into a monolithic system, we obtain
\begin{linenomath*}
\begin{equation}
\mathcal{C}_{CD} \mathbf{U}=
\begin{bmatrix}
    A_0 & M &  &   &  \\
    M & A_0 &M &   &  \\
     & M & A_0 & \ddots &  \\
     &  & \ddots & \ddots  & M\\
     &  &  &  M & A_0\\
\end{bmatrix}
\begin{bmatrix}
    \mathbf{u}_1  \\
    \mathbf{u}_2   \\
    \mathbf{u}_3   \\
    \vdots  \\
    \mathbf{u}_\ell
\end{bmatrix}
=
\begin{bmatrix}
	-M\mathbf{u}_0 \\
    \mathbf{0} \\
    \mathbf{0}    \\
    \vdots  \\
    \mathbf{0} 
\end{bmatrix} = \mathbf{b}_{CD},
\end{equation}
\end{linenomath*}

where $A_0 = \tau^2 K-2M$ for this problem. We can then precondition this system with the block circulant
\begin{linenomath*}
\begin{equation}
\mathcal{R}_{CD}=
\begin{bmatrix}
    A_0 & M &  &   &M  \\
    M & A_0 &M &   &  \\
     & M & A_0 & \ddots &  \\
     &  & \ddots & \ddots  & M\\
    M &  &  &  M & A_0\\
\end{bmatrix}.
\end{equation}
\end{linenomath*}

Alternatively we can approximate the second time derivative as
\begin{linenomath*}
\begin{equation}
\frac{\text{d}^2 \mathbf{u}_n}{\text{d} t^2} \approx \frac{\mathbf{u}_{n-2}-2\mathbf{u}_{n-1}+\mathbf{u}_n}{\tau^2}.
\label{BD2}
\end{equation}
\end{linenomath*}

Expression (\ref{BD2}) will be referred to as the 2-Step Backwards Difference (BD2) formula, which is a first order approximation.
By substituting (\ref{BD2}) into (\ref{applyEuler2}) and rearranging we obtain
\begin{linenomath*}
\begin{equation}
Mu_{n-2}-2Mu_{n-1}+(M + \tau^2K)u_n = \mathbf{0},
\end{equation}
\end{linenomath*}

which can be compiled into a monolithic system
\begin{linenomath*}
\begin{equation}
\mathcal{C}_{BD2} \mathbf{U}=
\begin{bmatrix}
    A_0 &  &  &   &  \\
    A_1 & A_0 & &   &  \\
    A_2 & A_1 & A_0 & &  \\
     & \ddots & \ddots & \ddots  & \\
     &  & A_2 &  A_1 & A_0\\
\end{bmatrix}
\begin{bmatrix}
    \mathbf{u}_1  \\
    \mathbf{u}_2   \\
    \mathbf{u}_3   \\
    \vdots  \\
    \mathbf{u}_\ell
\end{bmatrix}
=
\begin{bmatrix}
	(M+\tau^2K)\mathbf{u}_0 \\
    -M\mathbf{u}_0 \\
    \mathbf{0}    \\
    \vdots  \\
    \mathbf{0}
\end{bmatrix} = \mathbf{b}_{BD2},
\end{equation}
\end{linenomath*}

where $A_0 = M+\tau^2K, A_1 = -2M, A_2 = M$ for this problem. We will precondition this system with 
\begin{linenomath*}
\begin{equation}
\mathcal{R}_{BD2} =
\begin{bmatrix}
    A_0 &  &  & A_2  & A_1 \\
    A_1 & A_0 & &   & A_2 \\
    A_2 & A_1 & A_0 & &  \\
     & \ddots & \ddots & \ddots  & \\
     &  & A_2 &  A_1 & A_0\\
\end{bmatrix}.
\label{wavepreconditioner}
\end{equation}
\end{linenomath*}

We can also use a second order backwards difference formula
\begin{linenomath*}
\begin{equation}
\frac{\text{d}^2 \mathbf{u}_n}{\text{d} t^2} \approx \frac{-\mathbf{u}_{n-3}+4\mathbf{u}_{n-2}-5\mathbf{u}_{n-1}+2\mathbf{u}_n}{\tau^2},
\label{BD4}
\end{equation}
\end{linenomath*}

which is a second order method and will be referred to as the 4-Step Backwards Difference (BD4) formula. By substituting (\ref{BD4}) into (\ref{applyEuler2}) we obtain
\begin{linenomath*}
\begin{equation}
(2M+\tau^2K)\mathbf{u}_n - 5M\mathbf{u}_{n-1} + 4M\mathbf{u}_{n-2} -M\mathbf{u}_{n-3} =\mathbf{0}.
\end{equation}
\end{linenomath*}

As a consequence of using a 4-step method we have to approximate $u_1,u_2$ using sub-4-step methods. In this case we can use BD2 and the initial conditions. This results in the monolithic system
\begin{linenomath*}
\begin{equation}
\mathcal{C}_{BD4} = 
\begin{bmatrix}
   B &  &  & &  & &  & &  \\
    C&  B & & &  &  & & &  \\
    A_2 & A_1 & A_0 & & &  & & &  \\
    A_3 & A_2 & A_1 & A_0 & & & & & \\
     & A_3 & A_2 &  A_1 & A_0 & & & & \\
   & &\ddots & \ddots& \ddots & \ddots& & & \\
     & & & & & & & \\
     & & & &\ddots & \ddots& \ddots & \ddots&  \\
     &  & & & &   A_3 & A_2 &  A_1 & A_0\\
\end{bmatrix}
\begin{bmatrix}
    \mathbf{u}_1  \\
    \mathbf{u}_2   \\
    \mathbf{u}_3   \\
    \mathbf{u}_4   \\
    \mathbf{u}_5   \\
    \vdots  \\
       \\
\vdots  \\
    \mathbf{u}_\ell
\end{bmatrix}
=
\begin{bmatrix}
	B\mathbf{u}_0 \\
    -M\mathbf{u}_0 \\
    M\mathbf{u}_0     \\
    \mathbf{0}    \\
    \mathbf{0}    \\
    \vdots  \\
    \\
    \vdots  \\ 
    \mathbf{0}
\end{bmatrix} = \mathbf{b}_{BD4},
\end{equation}
\end{linenomath*}

where $B = M+\tau^2K, C = -2M, A_0 = 2M+\tau^2K, A_1 = -5M, A_2=4M, A_3=-M$ for this problem. 
We precondition this formulation of the wave equation system in a slightly different way. Previously we preconditioned the Toeplitz system with its corresponding circulant. In this instance, we do not have a Toeplitz structure to begin with. To overcome this, we precondition the system with the circulant matrix that would result if we had a Toeplitz system. That is, we precondition the system with the circulant matrix
\begin{linenomath*}
\begin{equation}
\mathcal{R}_{BD4} = 
\begin{bmatrix}
    A_0 &  &  & &  & &A_3  &A_2 &A_1  \\
    A_1&  A_0 & & &  &  & &A_3&A_2  \\
    A_2 & A_1 & A_0 & & &  & & &A_3  \\
    A_3 & A_2 & A_1 & A_0 & & & & & \\
     & A_3 & A_2 &  A_1 & A_0 & & & & \\
   & &\ddots & \ddots& \ddots & \ddots& & & \\
     & & & & & & & \\
     & & & &\ddots & \ddots& \ddots & \ddots&  \\
     &  & & & &   A_3 & A_2 &  A_1 & A_0\\
\end{bmatrix}.
\end{equation}
\end{linenomath*}

We will not pursue non-uniform temporal domains since there is less
interest in using non-uniform time-steps for the wave equation.

\section{Parallel Implementation}
\label{sec:parallel}
\smallskip

Throughout our implementation we keep all matrices on the master process and broadcast them when necessary. Vectors, on the other hand, are defined on all processes. To understand the reasoning for this, consider the fact that our dense block $U$ requires $\mathcal{O}(\ell^2)$ memory and our sparse blocks (linear combinations of $M,K$) each require $\mathcal{O}(n)$ memory. Since there are $\ell$ sparse blocks, the total memory cost is $\mathcal{O}(\ell^2 + n\ell)$. Further, since a vector requires $\mathcal{O}(n\ell)$ memory, if each one of $p$ processes has a copy of the vector, then we are using $O(n\ell p)$ memory. In the complexity analysis below, we consider the situation where $\ell$ processes are available. Using $\ell$ processes would significantly increase the memory requirements of our implementation with this memory management scheme. In our case, $p<<\ell$ and so the vector storage cost almost matches that of the matrix storage cost. Even in the case where $p \sim \ell$, the reduction in communication cost that results from this type of memory management is likely to be significant.  

In order to see how (\ref{standardheatprecon}) can be applied in parallel we follow \cite{McDonald_Pestana_W18} in writing it in the form
\begin{linenomath*}
\begin{equation}
\mathcal{P} = \mathbb{I}_{\ell} \otimes A_0 + \Sigma \otimes A_1,
\label{preconkron}
\end{equation}
\end{linenomath*}

where 
\begin{linenomath*}
\begin{equation}
\Sigma=
\begin{bmatrix}
   &   &        &   & 1\\
 1 &   &        &   &  \\
   & 1 &        &   &  \\
   &   & \ddots &   &   \\
   &   &        & 1 &
\end{bmatrix} \in \mathbb{R}^{n\times n},
\end{equation}
\end{linenomath*}

and $\mathbb{I}_\ell$ is the $\ell \times \ell$ identity matrix. The key property of circulant matrices is that they can be diagonalised by a Fourier basis. That is, we can write $\Sigma = U\Lambda U^*$, where $U_{k,j} = e^{((k-1)(j-1)\pi i )/n }/\sqrt{n}$ and the diagonal entries of $\Lambda$ are the roots of unity for the ``downshift" matrix $\Sigma$. Hence, again following \cite{McDonald_Pestana_W18}, (\ref{preconkron}) can be written as
\begin{linenomath*}
\begin{equation}
\mathcal{P} = (U \otimes \mathbb{I}_n)[\mathbb{I}_\ell \otimes A_0 + \Lambda \otimes A_1](U^{*} \otimes \mathbb{I}_n).
\end{equation}
\end{linenomath*}

Inverting $\mathcal{P}$ we obtain
\begin{linenomath*}
\begin{equation}
\mathcal{P}^{-1} = (U \otimes \mathbb{I}_n)[\mathbb{I}_\ell \otimes A_0 + \Lambda \otimes A_1]^{-1}(U^{*}\otimes \mathbb{I}_n).
\label{precon2para}
\end{equation}
\end{linenomath*}

The application of $(U\otimes \mathbb{I}_n)$ to a vector can be carried out in parallel. To see this, consider the explicit representation of $(U\otimes \mathbb{I}_n)$ given by
\begin{linenomath*}
\begin{equation}
(U\otimes \mathbb{I}_n) = 
\left[
\begin{array}{c|c|c}
U_{11}\mathbb{I}_n & \dots  &U_{1\ell}\mathbb{I}_n \\
\hline
\vdots & \ddots  &\vdots \\
\hline
U_{\ell 1}\mathbb{I}_n & \dots &  U_{\ell \ell}\mathbb{I}_n
\end{array}
\right].
\end{equation}
\end{linenomath*}

First, we broadcast each row of $U$ to a process. Then we can evaluate
\begin{linenomath*}
\begin{equation}
\mathbf{y}_i =
\left[
\begin{array}{c|c|c}
U_{i1}\mathbb{I}_n &  \dots &U_{i\ell}\mathbb{I}_n
\end{array}
\right]
\begin{bmatrix}
    \mathbf{z}_1  \\
    \vdots \\
    \mathbf{z}_\ell
\end{bmatrix}
\label{intermediatecal}
\end{equation}
\end{linenomath*}

on each process, where $\mathbf{z}_i \in \mathbb{R}^n$ is a chunk of the $n \ell$-vector $\mathbf{z}$ and $i$ is an integer such that $i \in [1,\ell]$. This can be done in $\mathcal{O}(n\ell)$. The local resultants of each of the calculations carried out on each process $\mathbf{y}_i$ can then be reduced to yield the final resultant of the calculation $(U\otimes \mathbb{I}_n)\mathbf{z}$. 

The only other implementation-specific issue that we need to be concerned with is the inversion of the block-tridiagonal matrix $\mathbb{I}_\ell \otimes A_0 + \Lambda \otimes A_1$. This can easily be carried out in parallel by assigning each tridiagonal block to a process and then applying the Thomas algorithm to each block, which would incur a cost of $\mathcal{O}(n^2)$ over $\ell$ processes. Therefore, the total complexity of this implementation is $\mathcal{O}(n^2+n\ell)$ over $\ell$ processes. 
As in \cite{McDonald_Pestana_W18}, multilevel iterations can be
applied as approximate solvers for the spatial operators when there 
is more than one spatial dimension.  

An alternative method involves the use of the Fast Fourier Transform (FFT). In order to see how the FFT can be applied in this case, we can write
\begin{linenomath*}
\begin{align}
(U\otimes \mathbb{I}_n) \mathbf{z}
&=
\left[
\begin{array}{c|c|c}
U_{11}\mathbb{I}_n & \dots  &U_{1\ell}\mathbb{I}_n \\
\hline
\vdots & \ddots  &\vdots \\
\hline
U_{\ell 1}\mathbb{I}_n & \dots &  U_{\ell \ell}\mathbb{I}_n
\end{array}
\right]
\begin{bmatrix}
    \mathbf{z}_1  \\
    \vdots \\
    \mathbf{z}_\ell
\end{bmatrix} \\
&=
\left[
\begin{array}{c|c|c}
U &   & \\
\hline
 & \ddots  &\\
\hline
 &  &  U
\end{array}
\right]
\begin{bmatrix}
    \mathbf{x}_1  \\
    \vdots \\
    \mathbf{x}_n
\end{bmatrix} \\
&=
(\mathbb{I}_n \otimes U)\mathbf{x},
\end{align}
\end{linenomath*}

where $(\mathbf{x}_i)_k = (\mathbf{z}_k)_i$ and $\mathbf{x}_i \in \mathbb{R}^\ell$ are the chunks of $\mathbf{x}$, $i\in[1,n]$ and $k \in [1,\ell]$. This transformation is called a vector transpose and can be carried out in $\mathcal{O}(n)$ over $\ell$ processes.\footnote{See \texttt{VectorTranspose} of \texttt{ParallelRoutines.cpp} in the Git repository.} Therefore we can form
\begin{linenomath*}
\begin{equation}
\mathbf{y}_i =
U\mathbf{x}_i,
\label{alternativecal}
\end{equation}
\end{linenomath*}

where $i$ is an integer such that $i\in [1,n]$. According to \cite{VANLOAN}, the cost of applying $U$ to a vector using the FFT is $\mathcal{O}(\ell \log \ell)$. Before we can progress with the calculation, we will have to apply the vector transpose again to $[\mathbf{y}_1,...,\mathbf{y}_n]^T$. Consequently the complexity of this implementation now becomes $\mathcal{O}(\ell \log \ell + n^2)$ over $\max(n,\ell)$ processes.

The timing results in the next section were obtained using the first method of evaluating $(U \otimes \mathbb{I}_n)$, not the FFT method. The reason for this is that the extra communication required to perform the transpose operator multiple times significantly reduced the performance of the all-at-once implementation. However, the functionality to implement both routines is provided in the GitHub repository.

From an inspection of the operations we can see that the most
expensive component of a GMRES iteration is the application of the
preconditioner. If we consider applying the above ideology to
$\mathcal{Q}_i^{-1}\mathcal{B}\mathbf{U} =
\mathcal{Q}_i^{-1}\mathbf{b}$, then for us to increase $i$ from 1 to 2
we have to take into account the cost of one extra preconditioner
application, as well as communication. That is, for us to consider the
Neumann method effective,
we must expect the number of GMRES iterations that are needed to solve the problem to reduce by more than half.

\section{Numerical results}
\label{sec:results}
\smallskip

All parallel results in this section were generated on the nightcrawler workstation at Oxford University. This machine is equipped with $2 \times 18$ core Intel(R) Xeon(R) Gold 6140 CPU @ 2.30GHz processors, 768GB RAM and 4600GB in scratch disk capacity. Care was taken to access the machine when the workload was low so as to maintain consistent results.

In \cite{McDonald_Pestana_W18} we see that the idea of preconditioning a block Toeplitz system with a block circulant matrix yields very good theoretical results. Table 1 of \cite{McDonald_Pestana_W18} shows that few iterations are required when computing the solution of such preconditioned Toeplitz systems. It was suggested in \cite{McDonald_Pestana_W18} that the all-at-once method can be executed in parallel. Timing results for the heat equation on a uniform temporal domain with initial condition $u_0^{s1}(x) = x(1-x)$ are given in Table \ref{table1}.
\begin{table}[h]
\caption{Timed results (in seconds) for solving the system $\mathcal{P}^{-1}\mathcal{A}\mathbf{U} = \mathcal{P}^{-1}\mathbf{b}$ using GMRES with tolerance set to $10^{-5}$. The iteration count remained at a constant value of 2 for all values of $n$ and $\ell$ tested. $p$ is the number of processes used in the calculations.}
\begin{center}
\begin{tabular}{|c c|c|c|c|c|c|c|}\hline
          &		    & $ p = 1 $ & $ p = 2 $ &$ p = 4 $  &$ p = 8 $ &$ p = 16 $ & $ p = 32 $ \\
\hline
          & $n=320$ & 77.72& 29.26& 15.32&  8.95& 5.11& 3.34\\
 $\ell= 768$ & $n=512$ & 152.64& 57.54& 32.71& 17.52&   11.54& 6.69\\
          & $n=768$ & 245.47& 97.77& 50.81& 30.71& 16.66& 9.65\\
 \hline
          & $n=320$ & 146.67&54.68&28.40 & 17.059 & 10.35 & 6.07\\
 $\ell= 1024$ & $n=512$ & 265.22&107.07&60.86 & 34.13 & 20.40 & 11.75 \\
          & $n=768$ & 459.12&198.94&101.23 & 55.85 & 28.55 & 16.12 \\
  \hline
          & $n=320$ & 325.14 & 124.67& 63.64   &39.78 & 22.74 & 13.06\\
 $\ell= 1440$& $n=512$ & 646.81 & 239.65 & 123.44 &72.44 & 40.95& 22.50 \\
          & $n=768$ & 979.85 & 432.46 & 215.77 &114.99 & 59.80 & 32.41\\
            \hline
 $\ell= 1440$& $n=1568$ & 2119.91&815.93&431.13&218.24&118.62&63.30\\
\hline
\end{tabular}
\end{center}
\label{table1}
\end{table}
\smallskip
Referring to Table \ref{table1} we can see that increasing the number of processes from 1 to 32 results in a significant reduction in the time taken to solve the preconditioned system associated with the heat equation. The most significant speed-up is achieved when $\ell=1440$ and $n=1568$: The time taken to solve the system reduces from $\sim 35$ minutes to $\sim 1$ minute. Observe that, for all values of $n$ and $\ell$, the time taken to solve the problem reduces by more than half as a result of increasing $p$ from 1 to 2. We suspect that this is because half of the problem fits in local memory better than the entire problem does.

While the data presented in Table \ref{table1} is very useful for highlighting the speed-up achieved by distributing the calculation across multiple processes, it would be useful to see how efficiently the processes are being used. To see this, we consider the parallel efficiency of $p$ processes defined by
\begin{linenomath*}
\begin{equation}
P_{eff}^p = \frac{\text{Time Taken on 1 Process}}{p \times \text{Time Taken on } p \text{ Processes}}.
\end{equation}
\end{linenomath*}

\begin{figure}[h]
\caption{The parallel efficiency of our implementation of GMRES used to solve the all-at-once formulation of the preconditioned heat equation system $\mathcal{P}^{-1}\mathcal{A}\mathbf{U} = \mathcal{P}^{-1}\mathbf{b}$. }
\centering
\includegraphics[scale=0.8]{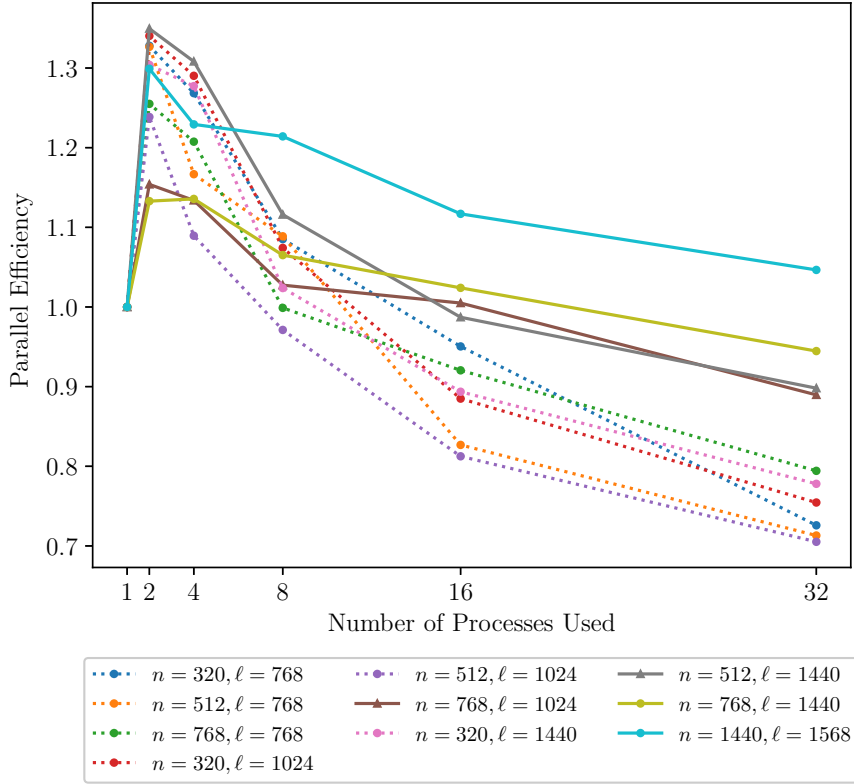}
\label{paraeff1}
\end{figure}

The parallel efficiency results are shown in Figure \ref{paraeff1}. The suspected reason for the jump in going from $P_{eff}^1$ to $P_{eff}^2$ is that, as noted above, the problem fits better in local memory over two processes. The parallel efficiency falls as we increase the number of processes used in the calculation, which is to be expected. This is a consequence of the number of communications taking place between processes. We note that as we increase the number of degrees of freedom, $P_{eff}^{32}$ also increases, particularly for the values $n=1568$ and $\ell=1440$, $P_{eff}^{32} > 1$. That is, our metric for measuring parallel efficiency implies that the all-at-once implementation is more efficient on 32 processes than it is on a single process. 

In Section 2 we introduced an extension to the all-at-once method that enables us to consider problems that are non-uniformly discretised in time. The key behind the extension is the Neumann approximation of the preconditioner $\mathcal{Q}^{-1}$, given by (\ref{preconapprox1}). The non-uniform temporal discretisation that was used to obtain the numerical results associated with the system $\mathcal{Q}_i^{-1}\mathcal{B}\mathbf{U} = \mathcal{Q}_i^{-1}\mathbf{b}$ is given by
\begin{linenomath*}
\begin{equation}
t_j = 
\begin{cases} 
	0, & j=0, \\
     \frac{1}{n}(j + \delta(\texttt{Rand(0,1,j)}-0.5)), & j \in [1,\ell-1], \\
      1, & j=\ell, 
   \end{cases}
\end{equation}
\end{linenomath*}

where $j$ is an integer, $\delta$ is a real number such that $\delta
\in (0,1)$ and \texttt{Rand(0,1,j)} is a random real number between
$0$ and $1$. Larger values of $\delta$ clearly tend to give more
irregular time-steps $\tau_j=t_j-t_{j-1}$. Table \ref{table2} shows the GMRES iteration counts required to solve $\mathcal{Q}_i^{-1}\mathcal{B}\mathbf{U} = \mathcal{Q}_i^{-1}\mathbf{b}$ for increasing values of $i$. We also show the difference between the time taken to solve the problem with $i=2$ and the time taken to solve the problem with $i=1$, $\Delta_{2,1}$. That is, a positive value of $\Delta_{2,1}$ indicates that it took longer to solve the problem with $i=2$ than it did with $i=1$. 
\begin{table}[h]
\caption{GMRES iteration count for $\mathcal{Q}_i^{-1}\mathcal{B}\mathbf{U} = \mathcal{Q}_i^{-1}\mathbf{b}$. The tolerance was set to $10^{-5}$. The values of $n$ and $\ell$ for each table are given as (a) $n=320,\ell=768$, (b) $n=512,\ell=768$, (c) $n=768,\ell=768$, (d) $n=512,\ell=1024$, (e) $n=768,\ell=1024$, (f) $n=1024,\ell=1024$. $\Delta_{1,2}$ is the difference between the time taken to solve the problem with $i=2$ and the time taken to solve the problem using $i=1$. In this particular experiment, 16 processes were utilised.}
\begin{center}
\begin{tabular}{|l|c|c|c|c||l|c|c|c|c|}\hline
(a)& $i=1$ & $i=2$ &$i=3$ & $\Delta_{2,1}$ & (b)& $i=1$ & $i=2$ &$i=3$ & $\Delta_{2,1}$\\
\hline
 $\delta = 0.9$ & 6 & 4 & 3 & 9.32 & $\delta = 0.9$ & 6 & 5 & 3 & 25.28 \\
 $\delta = 0.8$ & 6 & 4 & 2 & 9.36 &  $\delta = 0.8$ & 6 & 5 & 2 & 29.73\\
 $\delta = 0.7$ & 4 & 4 & 2 & 12.83 &  $\delta = 0.7$ & 4 & 4 & 2 & 20.83 \\
 $\delta = 0.6$ & 4 & 4 & 2 & 12.38 &  $\delta = 0.6$ & 4 & 3 & 2 & 21.38\\
 $\delta = 0.5$ & 4 & 4 & 2 & 11.75 &  $\delta = 0.5$ & 4 & 3 & 2 & 16.75 \\
 $\delta = 0.4$ & 4 & 4 & 2 & 10.33 &  $\delta = 0.4$ & 4 & 3 & 2 & 15.33 \\
 $\delta = 0.3$ & 4 & 3 & 2 & 7.90 &  $\delta = 0.3$ & 4 & 3 & 2 & 14.90\\
 $\delta = 0.2$ & 4 & 3 & 2 & 5.61 &  $\delta = 0.2$ & 4 & 3 & 2 & 15.61\\
 $\delta = 0.1$ & 4 & 3 & 2 & 4.32 &  $\delta = 0.1$ & 4 & 3 & 2 & 15.32\\
\hline
\end{tabular}
\begin{tabular}{|l|c|c|c|c||l|c|c|c|c|}\hline
(c)& $i=1$ & $i=2$ &$i=3$ & $\Delta_{2,1}$ & (d)& $i=1$ & $i=2$ &$i=3$ & $\Delta_{2,1}$\\
\hline
 $\delta = 0.9$ & 6 & 4 & 3 & 34.13 & $\delta = 0.9$ & 6 & 5 & 3 & 62.65 \\
 $\delta = 0.8$ & 6 & 4 & 2 & 33.75 &  $\delta = 0.8$ & 6 & 4 & 3 & 51.10\\
 $\delta = 0.7$ & 6 & 3 & 2 & 21.68 &  $\delta = 0.7$ & 4 & 4 & 3 & 59.12 \\
 $\delta = 0.6$ & 4 & 3 & 2 & 31.87 &  $\delta = 0.6$ & 4 & 4 & 3 & 58.01\\
 $\delta = 0.5$ & 4 & 3 & 2 & 33.43 &  $\delta = 0.5$ & 4 & 3 & 2 & 43.71 \\
 $\delta = 0.4$ & 4 & 3 & 2 & 32.42 &  $\delta = 0.4$ & 4 & 3 & 2 & 42.33 \\
 $\delta = 0.3$ & 4 & 3 & 2 & 31.31 &  $\delta = 0.3$ & 4 & 3 & 2 & 43.51\\
 $\delta = 0.2$ & 4 & 3 & 2 & 31.91 &  $\delta = 0.2$ & 4 & 3 & 2 & 42.51\\
 $\delta = 0.1$ & 4 & 3 & 2 & 32.36 &  $\delta = 0.1$ & 4 & 3 & 2 & 44.32\\
\hline
\end{tabular}
\begin{tabular}{|l|c|c|c|c||l|c|c|c|c|}\hline
(e)& $i=1$ & $i=2$ &$i=3$ & $\Delta_{2,1}$ & (f)& $i=1$ & $i=2$ &$i=3$ & $\Delta_{2,1}$\\
\hline
 $\delta = 0.9$ & 6 & 4 & 3 & 59.11 & $\delta = 0.9$ & 6 & 4 & 3 & 82.11\\
 $\delta = 0.8$ & 6 & 4 & 3 & 58.12 &  $\delta = 0.8$ & 4 & 4 & 3 & 103.12\\
 $\delta = 0.7$ & 4 & 3 & 3 & 63.22 &  $\delta = 0.7$ & 4 & 4 & 2 & 102.18 \\
 $\delta = 0.6$ & 4 & 3 & 2 & 65.26 &  $\delta = 0.6$ & 4 & 3 & 2 & 80.31\\
 $\delta = 0.5$ & 4 & 3 & 2 & 63.31 &  $\delta = 0.5$ & 4 & 3 & 2 & 75.42 \\
 $\delta = 0.4$ & 4 & 3 & 2 & 64.23 &  $\delta = 0.4$ & 4 & 3 & 2 & 74.15 \\
 $\delta = 0.3$ & 4 & 3 & 2 & 67.53 &  $\delta = 0.3$ & 4 & 3 & 2 & 78.13\\
 $\delta = 0.2$ & 4 & 3 & 2 & 64.21&  $\delta = 0.2$ & 4 & 3 & 2 & 74.21\\
 $\delta = 0.1$ & 4 & 3 & 2 & 65.26 &  $\delta = 0.1$ & 4 & 3 & 2 & 75.26\\
\hline
\end{tabular}
\end{center}
\label{table2}
\end{table}
\smallskip
Firstly, $\Delta_{2,1}>0$ for all values of $n,\ell$ and $\delta$ tested. This somewhat validates the claim we made in the previous section according to which the iteration count would need to reduce by more than half for any improvements to be observed in the timed results. For $n = 768,\ell=768$ and $\delta = 0.7$, the iteration count reduced by a half as a result of increasing $i$ from 1 to 2; in this case $\Delta_{2,1}=21.68$, which is quite significant. On a positive note, Table \ref{table2} highlights the robustness of the all-at-once formulation in the presence of temporal perturbations. For instance, increasing $\delta$ from 0.1 to 0.9 increases the iteration count by 2 for all values of $n$ and $\ell$ tested. Increasing $i$ beyond 2 does halve the iteration count in some cases, but since the cost incurred by communication is going to be increased further, there will be no advantage in doing so.

In Section 2 we introduced an all-at-once formulation of the wave equation. Three formulas were considered as candidates to approximate the time derivative. The central differences formulation performed poorly by selecting the smooth initial condition $u_0^{s2} = \sin(2 \pi x)$ and GMRES failed to converge to a solution by selecting the non-smooth initial condition
\begin{linenomath*}
\begin{equation}
u_0^{ns}(x)=
\begin{cases} 
      \cos^2 4 \pi \big(x-\frac{1}{2}\big), & x \in \big(\frac{3}{8},\frac{5}{8}\big), \\
      0, & x \in [0,1] \backslash \big(\frac{3}{8},\frac{5}{8}\big), 
   \end{cases}
\end{equation}
\end{linenomath*}

when applied to the system $\mathcal{R}_{CD}^{-1}\mathcal{C}_{CD} \mathbf{U} = \mathcal{R}_{CD}^{-1} \mathbf{b}_{CD}$. Central differences performed poorly on $u_0^{s2}$ in the sense that the solution displays aggressive numerical dissipation, regardless of how large $n$ and $\ell$ are chosen to be. However, the number of GMRES iterations required to solve the system $\mathcal{R}_{CD}^{-1}\mathcal{C}_{CD} \mathbf{U} = \mathcal{R}_{CD}^{-1} \mathbf{b}_{CD}$ with the smooth initial condition remained at 2 for all values of $n$ and $\ell$ tested. For these reasons the central differences formulation will not be considered further.

\begin{figure}[h]
\caption{Solution profiles for the wave equation. The profiles were taken at equally spaced time intervals and the non-smooth initial condition $u_0^{ns}$ was used. $n=\ell=128$. (a) $\mathcal{R}_{BD2}^{-1}\mathcal{C}_{BD2} \mathbf{U} = \mathcal{R}_{BD2}^{-1}\mathbf{b}_{BD2}$ (b) $\mathcal{R}_{BD4}^{-1}\mathcal{C}_{BD4} \mathbf{U} = \mathcal{R}_{BD4}^{-1}\mathbf{b}_{BD4}$. }
\centering
\includegraphics[scale=0.35]{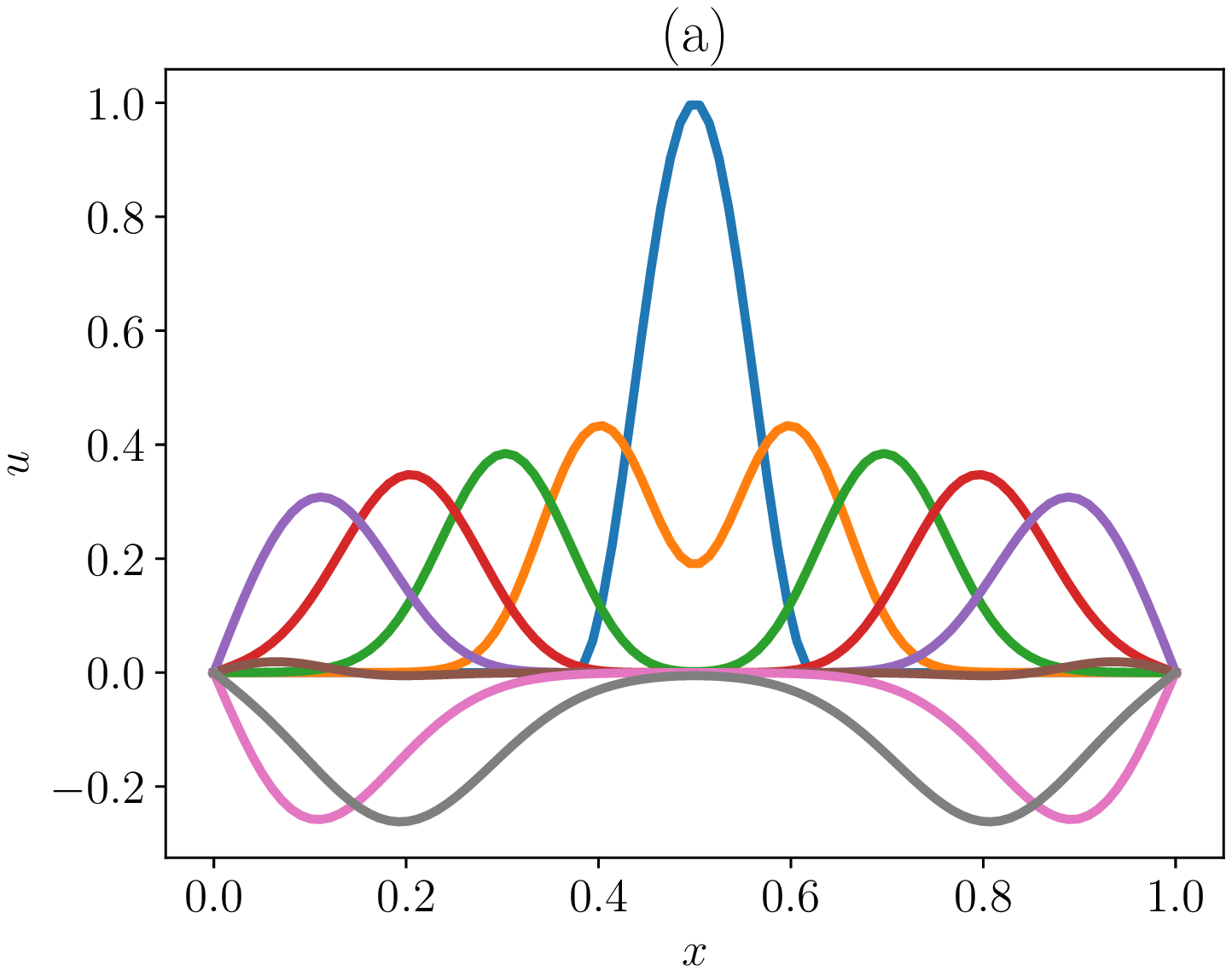}
\includegraphics[scale=0.35]{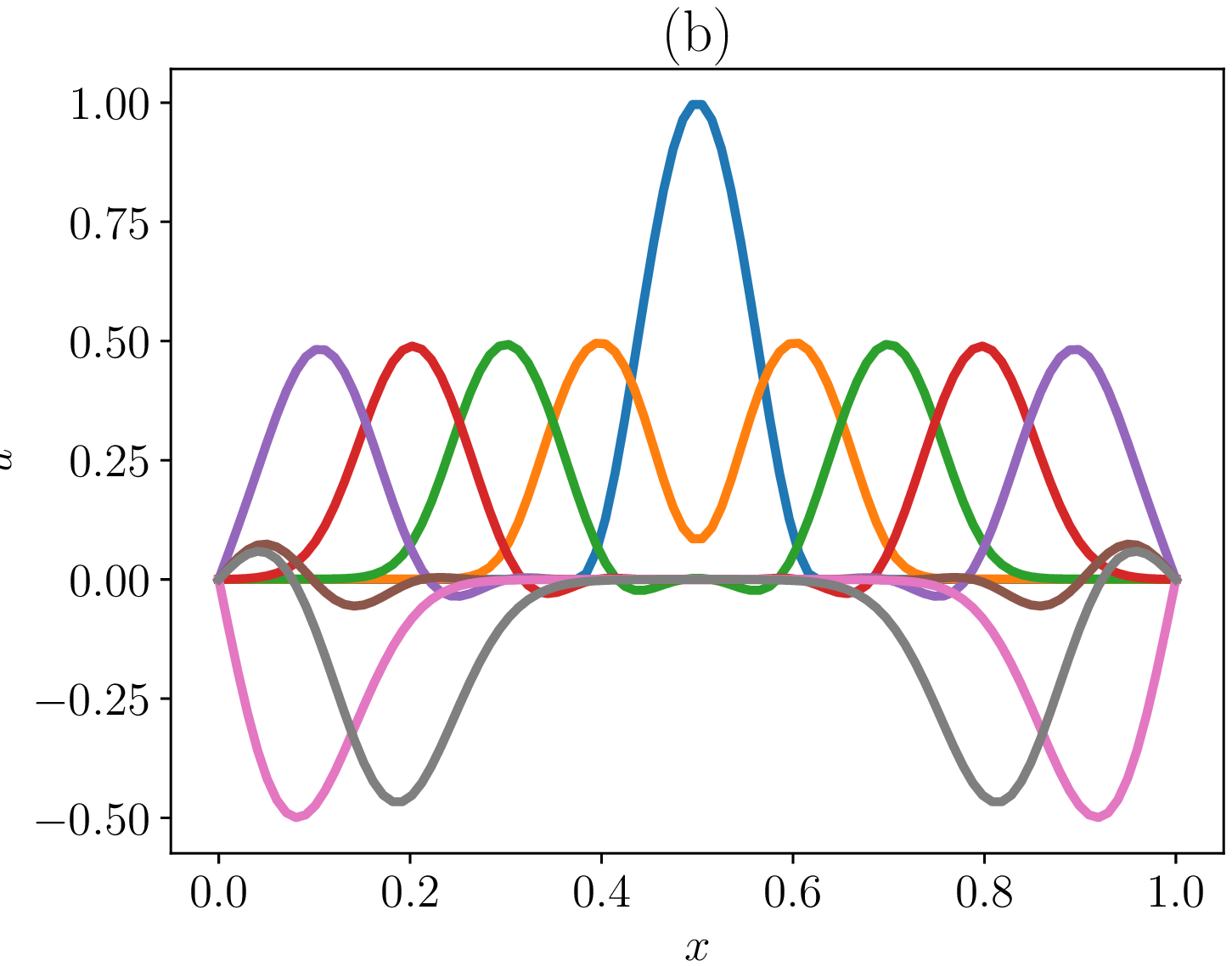}
\label{waves}
\end{figure}

Iteration counts for BD2 and BD4 are shown in Table \ref{table3} and Table \ref{table4}. Plots of the solution of the wave equation, obtained using the corresponding all-at-once formulations, are shown in Figure \ref{waves}. All of these results were obtained with the non-smooth initial condition $u_0^{ns}$. Iteration counts for BD2 are low and the wave speeds are largely conserved. The drawback of BD2 is the introduction of numerical dissipation in the solution, as seen in Figure \ref{waves} (a), although the dissipation is very subtle compared to that observed in the solution of the central differences formulation. BD4 rectifies this drawback, but iteration counts are much higher. Similarly to the central difference approximation, when the smooth initial condition was used, the number of GMRES iterations required to solve the problem remained fixed at 2 for both BD2 and BD4 formulations. Our observations indicate that the all-at-once method formulation of the wave equation is sensitive to the choice of initial conditions.

\begin{table}[h]
\caption{GMRES iteration counts $k$ required to solve the wave equation system. $v$ is the wave speed of the approximate solution. The wave speed of the exact solution is 1.
(a) $\mathcal{R}_{BD2}^{-1}\mathcal{C}_{BD2} \mathbf{U} = \mathcal{R}_{BD2}^{-1}\mathbf{b}_{BD2}$, (b) $\mathcal{R}_{BD4}^{-1}\mathcal{C}_{BD4} \mathbf{U} = \mathcal{R}_{BD4}^{-1}\mathbf{b}_{BD4}$. Tolerance was set to $10^{-5}$.}
\begin{center}
\begin{tabular}{|c c|c|c||c c|c|c|}\hline
      (a)    &		    &$k$ &  $v$  &(b) & &$k$ &  $v$  \\
\hline
          & $n=32$ & 5 & 1.03 &          &$n=32$ & 60 & 1.03\\
 $\ell= 32 $ & $n=64$ & 5 & 1.01 & $\ell= 32 $ & $n=64$& 100 & 1.01 \\
          & $n=96$ & 5 & 1.01 &          &  $n=96$ & 180& 0.67\\
\hline
          & $n=32$ & 6 & 2.06 &          &$n=32$ & 80 & 1.03\\
 $\ell= 64 $ & $n=64$ & 6 & 1.01 & $\ell= 64 $ & $n=64$& 120 & 1.01 \\
          & $n=96$ & 6 & 1.34 &          &  $n=96$ & 200& 0.67\\
\hline
          & $n=32$ & 6 & 3.09 &          &$n=32$ & 140 & 3.10\\
 $\ell= 96 $ & $n=64$ & 8 & 1.52 & $\ell= 96 $ & $n=64$& 180 & 1.52 \\
          & $n=96$ & 6 & 1.01 &          &  $n=96$ & 9216& 1.01\\
\hline
\end{tabular}
\end{center}
\label{table3}
\end{table}

\begin{table}[h]
\caption{GMRES iteration counts $k$ required to solve $\mathcal{R}_{BD2}^{-1}\mathcal{C}_{BD2} \mathbf{U} = \mathcal{R}_{BD2}^{-1}\mathbf{b}_{BD2}$ for larger values of $n$ and $\ell$.  Tolerance was set to $10^{-5}$.}
\begin{center}
\begin{tabular}{|c c|c|}\hline
          &		    &$k$ \\
\hline
          & $n=320$ & 8 \\
 $\ell= 768 $ & $n=512$ & 8  \\
          & $n=768$ & 8 \\
\hline
          & $n=320$ & 8\\
 $\ell=1024  $ & $n=512$ & 8 \\
          & $n=768$ & 8 \\
\hline
          & $n=320$ & 8 \\
 $\ell= 1440 $ & $n=512$ & 8\\
          & $n=768$ & 8 \\
\hline
\end{tabular}
\end{center}
\label{table4}
\end{table}

\begin{figure}[h]
\caption{The parallel efficiency of our implementation of GMRES used to solve the all-at-once formulation of the preconditioned wave equation system $\mathcal{R}_{BD2}^{-1}\mathcal{C}_{BD2} \mathbf{U} = \mathcal{R}_{BD2}^{-1}\mathbf{b}_{BD2}$. The smooth initial condition $u_0^{s2}$ was used.}
\centering
\includegraphics[scale=0.8]{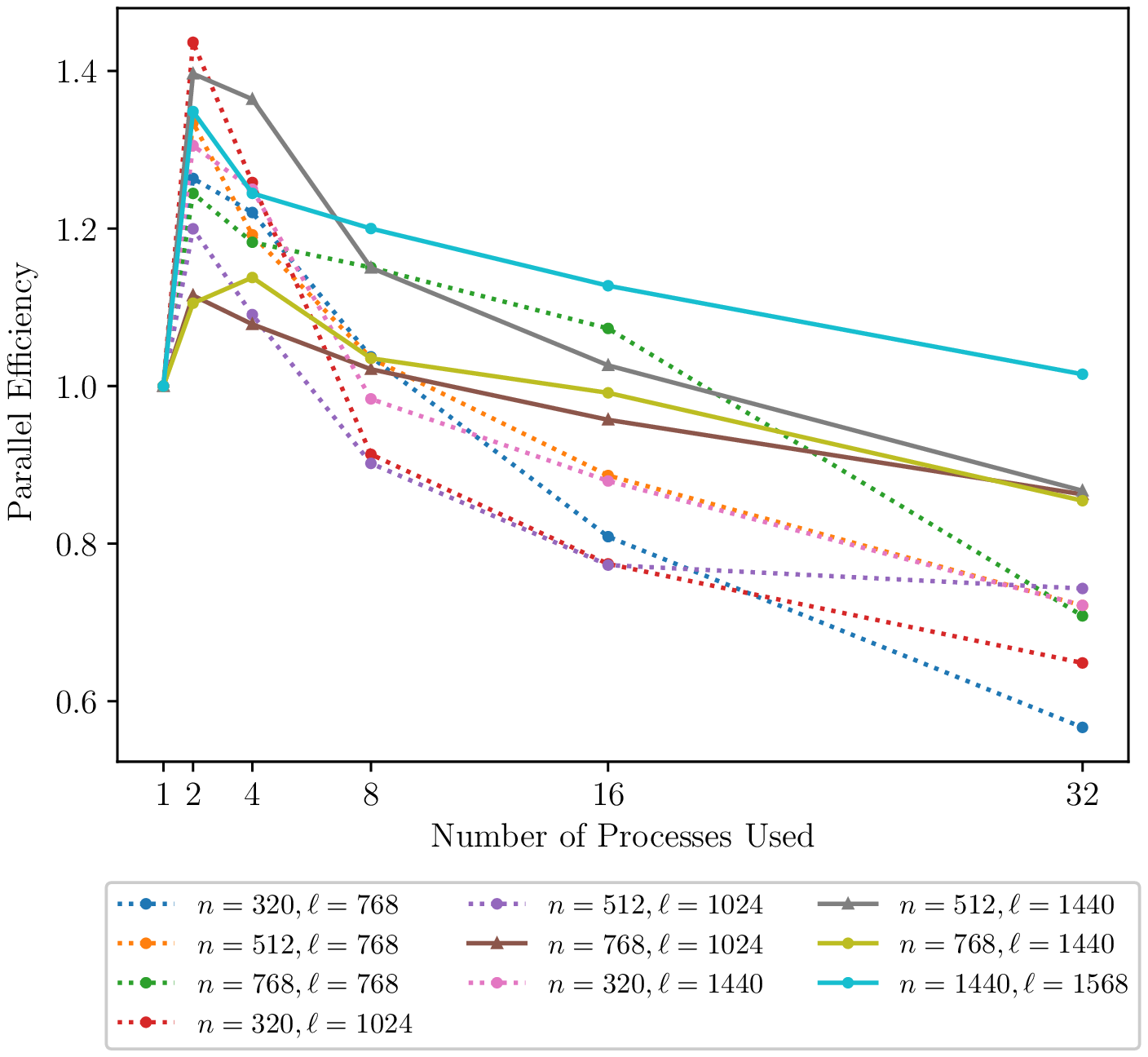}
\label{paraeff2}
\end{figure}

Timed results for the solution of the wave equation using the all-at-once formulation are provided in Table \ref{table5} and parallel efficiency results are displayed in Figure \ref{paraeff2}. The system associated with the wave equation resulting from using BD2 is less sparse than the system associated with the heat equation and so we would expect the computations to take longer. Also note that we have a reduction in parallel efficiency for all values of $n$ and $\ell$ tested. However, the parallel efficiency does seem to increase with the number of degrees of freedom. This trend was observed in the parallel efficiency results for the heat equation.

\begin{table}[h]
\caption{The timed results (in seconds) for solving the system $\mathcal{R}_{BD2}^{-1}\mathcal{C}_{BD2} \mathbf{U} = \mathcal{R}_{BD2}^{-1}\mathbf{b}_{BD2}$ using GMRES with tolerance set to $10^{-5}$. The iteration count remained at a constant value of 2 for all values of $n$ and $\ell$ tested. $p$ is the number of processes used in the calculations. The smooth initial condition $u_0^{s2}$ was used.}
\begin{center}
\begin{tabular}{|c c|c|c|c|c|c|c|}\hline
          &		    & $ p = 1 $ & $ p = 2 $ &$ p = 4 $  &$ p = 8 $ &$ p = 16 $ & $ p = 32 $ \\
\hline
          & $n=320$ & 79.07& 31.29& 16.20&  9.53& 6.11& 4.36\\
 $\ell= 768$ & $n=512$ & 163.68& 61.33& 34.33& 19.76&   11.54& 7.09 \\
          & $n=768$ & 251.37& 100.99& 53.14& 27.31& 14.64& 11.09\\
 \hline
          & $n=320$ & 153.37&53.39 &30.47 &20.99& 12.38&7.39\\
 $\ell= 1024$ & $n=512$ & 287.23&119.72&65.84 &39.81& 23.23&12.08\\
          & $n=768$ & 497.12&222.93&115.24&60.84& 32.46&18.01 \\
  \hline
          & $n=320$ &328.17&125.71&65.64  &41.70   &23.32&14.21\\
 $\ell= 1440$& $n=512$ &680.15&243.53&124.65 &73.92  &41.42&24.51 \\
          & $n=768$ &960.33&434.46&211.01 &115.95 &60.54&35.12\\
            \hline
 $\ell= 1440$& $n=1568$ & 2211.97 &820.21 &444.31 &230.42&122.63 &68.10\\
\hline
\end{tabular}
\end{center}
\label{table5}
\end{table}

\section{Conclusions}

We have provided a parallel implementation of the all-at-once method
of \cite{McDonald_Pestana_W18}. This was achieved using MPI and C++
and is a proof-of-concept software that supports the claims made in
\cite{McDonald_Pestana_W18}, namely, that the all-at-once method with
the McDonald et.al. preconditioner is parallelisable. We have also
provided new applications for the all-at-once method, namely,
applications to non-uniform temporal discretisation and to hyperbolic
equations. Problems in one spatial dimension were considered
throughout. To apply the all-at-once method to higher-dimensional
problems, some alterations to the implementation provided in this
paper are required. An alterative suggested in
\cite{McDonald_Pestana_W18} 
would be to apply Algebraic Multigrid in parallel instead of the 
parallel Thomas algorithm.

\vspace*{0.3in}


\end{document}